\documentclass[]{elsarticle}
\usepackage{amsmath,amssymb,amsfonts}
\usepackage{color}

\newtheorem{Th}{Theorem}[section]
\newtheorem{Le}{Lemma}[section]
\newtheorem{De}{Definition}[section]
\newtheorem{Rem}{Remark}[section]
\newtheorem{Cor}{Corollary}[section]

\newcommand{\bth}{\begin{Th}}
\newcommand{\eeth}{\end{Th}}
\newcommand{\ble}{\begin{Le}}
\newcommand{\eele}{\end{Le}}
\newcommand{\bde}{\begin{De}}
\newcommand{\ede}{\end{De}}
\newcommand{\bre}{\begin{Rem}}
\newcommand{\eere}{\end{Rem}}

\newcommand{\bco}{\begin{Cor}}
\newcommand{\eeco}{\end{Cor}}

\newcommand{\p}{\rm{Proof:$\,$}}

\newcommand{\beq}{\begin{equation}}
\newcommand{\eeq}{\end{equation}}

\begin{document}
\begin{frontmatter}

\title{Schwarz boundary value problems for polyanalytic equation in a sector ring}
\author[rvt]{Zhihua Du}
\author[rvtt]{ Ying Wang\corref{cor11}}
\cortext[cor11]{Corresponding author.\\ E-mail: wangyingyezi@sina.com.}
\author[rvttt]{Min Ku}

\address[rvt]{Department of Mathematics, Jinan University, Guangzhou 510632, China}
\address[rvtt]{School of Statistics and Mathematics, Zhongnan University of Economics and Law, Wuhan 430072, China}
\address[rvttt]{Department of Computing Science, University of Radboud, 6525 EC Nijmegen, Netherlands}

\begin{abstract} In this article, we
first give a modified Schwarz-Pompeiu formula in a general sector
ring with angle $\theta=\frac{\pi}{\alpha},\ \alpha\geq 1/2$ by
proper conformal mappings, and obtain the solution of the Schwarz problem for the Cauchy-Riemann equation in explicit forms. Furthermore, a class of integral operators is introduced together with their properties. Finally, by virtue of these operators, Schwarz problems for a inhomogeneous polyanalytic equation  and for a generalized polyanalytic equation are investigated, respectively.
\end{abstract}


\end{frontmatter}

\section{Introduction}
The theory of the classical boundary value problems for analytic functions has been applied directly or indirectly in many different fields [1-3], such as signal analysis, crack and elasticity, orthogonal polynomials, time-frequency and so on. This makes a great interest in the investigation of boundary value
problems for complex partial differential equations in different
domains [4-20]. Especially, in \cite{2020,V2008}, four basic boundary value problems for the Cauchy-Riemann equation and the Neumann problem for Bitsadze equation were studied in a ring domain by constructing kernel functions and integral expression formulas, respectively. In \cite{B2009}, some harmonic boundary value problems for the Poisson equation were investigated in upper half ring domain with the help of building a modified Cauchy-Pompeiu formula, harmonic Green and Neumann functions. Dirichlet and Neumann problems for the Poisson equation in a quarter sector ring and the general sector ring were respectively solved in \cite{W2015,B2012} through establishing proper Poisson kernels on basis of reflection principle.

 Motivated by these, our aim is to set up the theory of boundary value problems for polyanalytic equations and generalized polyanalytic equations in general sector rings through developing proper kernel functions, which generalizes the boundary value problem not only to the case of high order complex partial differential equation, but to more general sector ring domain. In this present paper, we first give a Schwarz-Pompeiu formula in the general sector ring with angle
$\theta=\frac{\pi}{\alpha}$, $\alpha\geq 1/2$ by proper conformal mappings, and obtain the solution of Schwarz problem for the Cauchy-Riemann equation. It is just the case in \cite{B2009} for $\alpha=1$ and the case in \cite{B2012} for $\alpha=2$. In Section 3, a poly-Schwarz and a poly-Pompeiu operator are introduced together with their properties, and then a Schwarz problem for the inhomogeneous polyanalytic equation in the general sector ring is investigated explicitly. In the end, the poly-Pompeiu operators are discussed in more detail for the upper half ring. On account of these poly-Pompeiu operators, we consider the Schwarz problem for the generalized inhomogeneous polyanalytic equations by transforming it into a singular integral equation.

Throughout this article, $\Omega$ is a sector ring with $\theta=\frac{\pi}{\alpha}$ $\ (\alpha\geq
1/2)$, that is, $\Omega=\{
 0<r<|z|<1,\ 0<\arg z<\frac{\pi}{\alpha},\ \alpha\geq 1/2\}$, whose boundary
 $\partial \Omega=[r,1]\cup \Gamma_1\cup [\varpi,\omega]\cup\Gamma_2$ is given counter-clockwisely, where
 four corner points are $r,\ 1,\ \varpi=e^{i\theta},\ \omega=re^{i\theta}$, respectively, and
$$ \Gamma_1:\ \tau\longmapsto e^{i \tau},\ \ \tau\in\left[0,\ \frac{\pi}{\alpha}\right],$$
$$ \Gamma_2:\ \tau\longmapsto re^{i \tau},\ \ \tau\in\left[\frac{\pi}{\alpha},\ 0\right].$$

\section{Schwarz Problem for the Cauchy-Riemann equation}
In order to discuss the Schwarz problem for the Cauchy-Riemann equation in $\Omega$, we define kernel functions $H_1(z,\zeta)$, $H_2(z,\zeta)$ as follows.
\begin{equation}\label{5.1}
H_1(z,\zeta)=\displaystyle\frac{\zeta+z}{\zeta-z}-\displaystyle\frac{\overline{\zeta}+z}{\overline{\zeta}-z}+2\sum\limits_{n=1}^\infty
r^{2\alpha n}\left[\displaystyle\frac{\zeta}{r^{2\alpha
n}\zeta-z}-\displaystyle\frac{z}{r^{2\alpha
n}z-\zeta}-\displaystyle\frac{\overline{\zeta}}{r^{2\alpha
n}\overline{\zeta}-z} +\displaystyle\frac{z}{r^{2\alpha
n}z-\overline{\zeta}}\right],\end{equation}
\begin{equation}\label{5.2}\begin{array}{ll}
H_2(z,\zeta)&=\displaystyle\frac{1}{\zeta-z}-\displaystyle\frac{z}{1-z\zeta}\\[3mm]&+\displaystyle\sum\limits_{n=1}^\infty
r^{2\alpha n}\left[\displaystyle\frac{1}{r^{2\alpha
n}\zeta-z}-\displaystyle\frac{z}{\zeta(r^{2\alpha
n}z-\zeta)}-\displaystyle\frac{1}{\zeta(r^{2\alpha n}-z\zeta)}
+\displaystyle\frac{z}{r^{2\alpha n}z\zeta-1}\right].\end{array}\end{equation}

\ble {\rm(\cite{B2009})}
Any \ $w\in C^1(R^+;\mathbb{C}
)\cap C(\overline{ R^+};\mathbb{C})$ can be expressed as
$$\begin{array}{ll} w(z)&=\displaystyle\frac{1}{2\pi
i}\int_{L_1\cup L_2}\mbox{\rm Re}
w(\zeta)H_1(z,\zeta)\frac{\rm d\zeta}{\zeta}+\frac{1}{\pi}
\int_{L_1}\frac{\mbox{\rm Im} w(\zeta)}{\zeta}\mbox{\rm
d}\zeta\\[3mm]&+\displaystyle\frac{1}{\pi
i}\int_{[-1,-r]\cup[r,1]}\mbox{\rm Re}
w(\zeta)H_2(z,\zeta)\mbox{\rm
d}\zeta-\displaystyle\frac{1}{\pi}\displaystyle\iint_{R^+}
\big[w_{\bar{\zeta}}(\zeta)H_2(z,\zeta)-\overline{w_{\overline{\zeta}}(\zeta)}H_2(z,\overline{\zeta})\big]\mbox{\rm
d}\xi \mbox{\rm d}\eta,
\end{array}
$$
 where $H_1(z,\zeta)$, $H_2(z,\zeta)$ are given by \eqref{5.1} and \eqref{5.2} with $\alpha=1$, \ $R^+=\{\,0<r<|z|<1, \ \mbox{\rm Im}z>0\}$, and its boundary $\partial R^+=[-1,-r]\cup [r,1]\cup L_1\cup L_2$ with $L_1=\{\,|\tau|=1,\ \mbox{\rm
 Im}\tau>0\}$ and
$L_2=\{\,|\tau|=r,\ \mbox{\rm
 Im}\tau>0\}$ being oriented counter-clockwise and clockwise, respectively.
\eele

By Lemma 2.1 and proper conformal mappings, we obtain the Schwarz-Pompeiu formula for the sector ring with angle
$\displaystyle\frac{\pi}{\alpha}\ (\alpha\geq 1/2)$ as follows.

 \bth Any \ $w\in C^1(\Omega;\mathbb{C}
)\cap C(\overline{ \Omega};\mathbb{C})$ has the following expression formula
\beq\label{3.1}
\begin{array}{ll}
w(z)=\displaystyle\frac{\alpha}{2\pi
i}\int_{\Gamma_1\cup\Gamma_2}\mbox{\rm Re}
w(\tau)H_1(z^\alpha,\tau^\alpha)\frac{\rm
d\tau}{\tau}+\displaystyle\frac{\alpha}{\pi
i}\int_{[\varpi,\omega]\cup[r,1]}\mbox{\rm Re}
w(\tau)H_2(z^\alpha,\tau^\alpha)\tau^{\alpha-1}\mbox{\rm
d}\tau\\[5mm]\ \ \ \ +\displaystyle\frac{\alpha}{\pi}
\int_{\Gamma_1}\frac{\mbox{\rm Im} w(\tau)}{\tau}\mbox{\rm
d}\tau
-\displaystyle\frac{\alpha}{\pi}\displaystyle\iint_{\Omega}
\big[\tau^{\alpha-1}H_2(z^\alpha,\tau^\alpha)f(\tau)-\overline{\tau^{\alpha-1}}H_2(z^\alpha,\overline{\tau^\alpha})\overline{f(\tau)}\
\big] \mbox{\rm d}\tau_1 \mbox{\rm d}\tau_2,
\end{array}
\eeq
where $\tau=\tau_1+i\ \tau_2,\ \tau_1,\ \tau_2\in\mathbb{R}$, $H_1$, $H_2$ are given by \eqref{5.1} and \eqref{5.2}, respectively.
\eeth

\p
Suppose $\widetilde{R^+}=\{\,0<r^\alpha<|z|<1, \ \mbox{\rm Im}z>0\}$, and
$\widetilde{L_2}=\{\,|\tau|=r^\alpha,\ \mbox{\rm
 Im}\tau>0\}$ is oriented clockwise. Consider a conformal mapping \cite{Wen1992},
$$\begin{array}{cc}
 \zeta: &\Omega\rightarrow \widetilde{R^+}\\
        &z\mapsto z^\alpha
        \end{array}$$
 where the branch is cut along $(-\infty,r)$, which maps $\Gamma_1$ onto
 $L_1$, $\Gamma_2$  onto
 $\widetilde{L_2}$, $[r,1]$ onto $[r^\alpha,1]$, and $[\varpi,\omega]$ is mapped onto $[-1,-r^\alpha]$, respectively. In the same time, the mapping
$$
 \begin{array}{cc}
 \varsigma: &\widetilde{ R^+}\rightarrow \Omega\\
        &z\mapsto z^{1/\alpha}
 \end{array}$$
transforms the boundary of  $\widetilde{R^+}$ onto the relative boundary of $\Omega$, respectively.

Let $\partial_{\overline{z}}w(z)=f(z)$, and define
$$
W(z)=w(z^{1/\alpha}),\ \ \ \ z\in \widetilde{R^+},$$ then, $W(z)\in
C^1(\widetilde{R^+};\mathbb{C} )\cap C(\overline{ \widetilde{R^+}};\mathbb{C})$.
According to Lemma 2.1, for $z\in \Omega$,
$$\begin{array}{ll}
W(z^\alpha)&=\displaystyle\frac{1}{2\pi i}\int_{L_1\cup
\widetilde{L_2}}\mbox{\rm Re} W(\zeta)H_1(z^\alpha,\zeta)\frac{\rm
d\zeta}{\zeta}+\displaystyle\frac{1}{\pi
i}\int_{[-1,-r^\alpha]\cup[r^\alpha,1]}\mbox{\rm Re}
W(\zeta)H_2(z^\alpha,\zeta)\mbox{\rm d}\zeta\\[7mm]&+\displaystyle\frac{1}{\pi}
\int_{L_1}\frac{\mbox{\rm Im} W(\zeta)}{\zeta}\mbox{\rm
d}\zeta-\displaystyle\frac{1}{\pi}\displaystyle\iint_{\widetilde{R^+}}
\Big[W_{\bar{\zeta}}(\zeta)H_2(z^\alpha,\zeta)-\overline{W_{\overline{\zeta}}(\zeta)}H_2(z^\alpha,\overline{\zeta})\Big]\mbox{\rm
d}\xi \mbox{\rm d}\eta,
\end{array}
$$
where
$W_{\bar{\zeta}}(\zeta)=\displaystyle\frac{1}{\alpha}\overline{\zeta^{\frac{1}{\alpha}-1}}\ f(\zeta^{1/\alpha}).$ Taking a transformation $\zeta=\tau^\alpha$ with $\tau=\tau_1+i \tau_2$, then $\mbox{\rm
d}\xi\mbox{\rm d}\eta=\alpha^2|\tau|^{2(\alpha-1)}\mbox{\rm
d}\tau_1\mbox{\rm d}\tau_2$. Therefore, when $z\in \Omega$,
$$
\begin{array}{ll}
w(z)&=W(z^\alpha)=\displaystyle\frac{\alpha}{2\pi
i}\int_{\Gamma_1\cup\Gamma_2}\mbox{\rm Re}
W(\tau^\alpha)H_1(z^\alpha,\tau^\alpha)\frac{\rm
d\tau}{\tau}+\frac{\alpha}{\pi}
\int_{\Gamma_1}\frac{\mbox{\rm Im} W(\tau^\alpha)}{\tau}\mbox{\rm d}\tau\\[5mm]
& \ \ \ +\displaystyle\frac{\alpha}{\pi
i}\int_{[\varpi,\omega]\cup[r,1]}\mbox{\rm Re}
W(\tau^\alpha)H_2(z^\alpha,\tau^\alpha)\tau^{\alpha-1}\mbox{\rm
d}\tau\\[5mm]
&\ \
-\displaystyle\frac{\alpha}{\pi}\displaystyle\iint_{\Omega}
\Big[\tau^{\alpha-1}H_2(z^\alpha,\tau^\alpha)f(\tau)-\overline{\tau^{\alpha-1}}H_2(z^\alpha,\overline{\tau^\alpha})\overline{f(\tau)}\Big]
\mbox{\rm d}\tau_1 \mbox{\rm d}\tau_2.
\end{array}
$$
Thus, the proof is completed.

From the Lemmas 2.1-2.8 in \cite{W2015}, we easily obtain the following Lemma.
\ble
When $\gamma(\zeta)\in C(\partial\Omega ;\mathbb{R})$ and
$t\in \partial\Omega$,
$$\lim\limits_{z\in \Omega,\ z\rightarrow t}\displaystyle\frac{\alpha}{2\pi i}\int_{\partial\Omega}\gamma(\zeta)
[H_2(z^\alpha,\zeta^\alpha)-H_2(\overline{z^\alpha},\zeta^\alpha)]\zeta^{\alpha-1}\mbox{\rm
d}\zeta=\gamma(t).$$
\eele

Combining Theorem 2.1 with Lemma 2.2, we can solve the following Schwarz problem.

 \bth The Schwarz problem for the Cauchy-Riemann equation in $\Omega$
$$\left\{\begin{array}{ll}
w_{\overline{z}}=f \ \text{in}\ \Omega,\\[3mm]
\mbox{\rm Re}w=\gamma \ \text{on}\ \partial \Omega,\\[3mm]
\displaystyle\frac{\alpha}{\pi i}\int_{\Gamma_1}\frac{\mbox{\rm
Im}w(\zeta)}{\zeta}\mbox{\rm d}\zeta=c,\ \ c\in \mathbb{R}\end{array}\right.$$
with
$f\in L_p(\Omega;\mathbb{C}),$ $p>2$ and $\gamma\in C(\partial
\Omega;\mathbb{R})$, is solvable by
\begin{equation}\label{11sec1}\begin{array}{ll}
w(z)&=\displaystyle\frac{\alpha}{2\pi
i}\int_{\Gamma_1\cup\Gamma_2}\gamma(\zeta)H_1(z^\alpha,\zeta^\alpha)\frac{\rm
d\zeta}{\zeta} +\displaystyle\frac{\alpha}{\pi
i}\int_{[\varpi,\omega]\cup[r,1]}\gamma(\zeta)H_2(z^\alpha,\zeta^\alpha)\zeta^{\alpha-1}\mbox{\rm d}\zeta\\[7mm]
&+i c-\displaystyle\frac{\alpha}{\pi}\displaystyle\iint_{\Omega}\Big[\zeta^{\alpha-1}H_2(z^\alpha,\zeta^\alpha)f(\zeta)
-\overline{\zeta^{\alpha-1}}H_2(z^\alpha,\overline{\zeta^\alpha})\overline{f(\zeta)}\Big]\mbox{\rm
d}\xi \mbox{\rm d}\eta,\ \ \ z\in \Omega,
\end{array}
\end{equation}
where $\alpha\geq \displaystyle\frac{1}{2}$,  and $H_1$, $H_2$ are defined by \eqref{5.1}, \eqref{5.2}, respectively.
\eeth

\p By Theorem 2.1, it is sufficient to prove that \eqref{11sec1}
is the solution of the Schwarz problem. Obviously, the boundary integral and constant in
\eqref{11sec1} are analytic. Moreover, the kernel in area integral can be rewritten
$$
\frac{\alpha \zeta^{\alpha-1}}{\zeta^\alpha-z^\alpha}+{\text
{others}}=\frac{1}{\zeta-z}+g(z,\zeta),\ \ \
 z,\ \zeta\in \Omega.$$
 Here,
$g(z,\zeta)$ is analytic about $z$. Thus, it is easy to obtain
$$\partial_{\overline{z}}w(z)=\partial_{\overline{z}}
\left\{-\displaystyle\frac{\alpha}{\pi}\displaystyle\iint_{\Omega}\frac{\zeta^{\alpha-1}f(\zeta)}{\zeta^\alpha-z^\alpha}\mbox{\rm
d}\xi \mbox{\rm d}\eta\right\}
=\partial_{\overline{z}}\left\{-\displaystyle\frac{1}{\pi}\displaystyle\iint_{\Omega}\displaystyle\frac{f(\zeta)}{\zeta-z}\mbox{\rm
d}\xi \mbox{\rm d}\eta\right\}=f(z).$$

When\ $\zeta\in \Gamma_1$,
$$\begin{array}{ll}
\displaystyle\frac{\alpha}{2\pi
i}\int_{\Gamma_1}H_1(z^\alpha,\zeta^\alpha)\frac{\mbox{\rm d}z}{z}\\[3mm]
=\displaystyle\frac{1}{2\pi
i}\int_{|z|=1}\Bigg(\frac{\zeta^\alpha+z}{\zeta^\alpha-z}
+2\sum_{m=1}^\infty\bigg[\frac{r^{2m\alpha}\zeta^\alpha}{r^{2m\alpha}\zeta^\alpha-z}
-\frac{r^{2m\alpha}\overline{\zeta^\alpha}}{r^{2m\alpha}\overline{\zeta^\alpha}-z}\bigg]\Bigg)\frac{\mbox{\rm
d}z}{z}=0.
\end{array}
$$
Further, for $\zeta\in \Gamma_2$, we have
$$\begin{array}{ll}
\displaystyle\frac{\alpha}{2\pi
i}\int_{\Gamma_1}H_1(z^\alpha,\zeta^\alpha)\frac{\mbox{\rm
d}z}{z}\\[3mm]
=\displaystyle\frac{1}{\pi
i}\int_{|z|=1}\Bigg(\frac{\zeta^\alpha}{\zeta^\alpha-z}
+\sum_{m=1}^\infty\bigg[\frac{r^{2m\alpha}\zeta^\alpha}{r^{2m\alpha}\zeta^\alpha-z}
-\frac{r^{2m\alpha}}{r^{2m\alpha}-z\zeta^\alpha}\bigg]\Bigg)\frac{\mbox{\rm
d}z}{z}=0.
\end{array}
$$
Also, for $\zeta\in(r,1)\cup(\varpi,\omega)\cup \Omega$,
$$\begin{array}{ll}
\displaystyle\frac{\alpha}{2\pi
i}\int_{\Gamma_1}H_2(z^\alpha,\zeta^\alpha)\frac{\mbox{\rm
d}z}{z}\\[3mm]
=\displaystyle\frac{1}{2\pi
i}\int_{|z|=1}\Bigg(\frac{1}{\zeta^\alpha-z}
+\sum_{m=1}^\infty\bigg[\frac{r^{2m\alpha}}{r^{2m\alpha}\zeta^\alpha-z}
-\frac{z}{\zeta^\alpha(z-r^{-2\alpha
m}\zeta^\alpha)}\bigg]\Bigg)\frac{\mbox{\rm
d}z}{z}=0.
\end{array}
$$
Similarly, when\ $\zeta\in \Omega$,
$$
\displaystyle\frac{\alpha}{2\pi
i}\int_{\Gamma_1}H_2(z^\alpha,\overline{\zeta^\alpha})\frac{\mbox{\rm
d}z}{z}=0.
$$
Therefore, by interchange of integral order, we obtain  \ $\displaystyle\frac{\alpha}{\pi
i}\int_{\Gamma_1}\frac{\mbox{\rm Im}w(\zeta)}{\zeta}\mbox{\rm
d}\zeta=c.$

Suppose $w_0$ be the area integral in \eqref{11sec1}, by simple calculation,
$$
\mbox{\rm
Re}w_0(z)=-\displaystyle\frac{\alpha}{2\pi}\iint_{\Omega}\left[f(\zeta)G^*(z,\zeta)+\overline{f(\zeta)}\
\overline{G^*(z,\zeta)}\right]\mbox{\rm d}\xi \mbox{\rm d}\eta$$
with
$G^*(z,\zeta)=\zeta^{\alpha-1}[H_2(z^\alpha,\zeta^\alpha)-H_2(\overline{z^\alpha},\zeta^\alpha)]$.
Since for $(z,\zeta)\in
\partial\Omega\times\Omega$, $G^*(z,\zeta)=0$, which means that $\mbox{\rm Re}w_0(z)=0 $ is true for $z\in
\partial\Omega$. Then, we rewrite
\beq\label{3.52}\begin{array}{ll}
\mbox{\rm Re}w(z)=\displaystyle\frac{\alpha}{2\pi
i}\int_{\partial\Omega}\gamma(\zeta)[H_2(z^\alpha,\zeta^\alpha)-H_2(\overline{z^\alpha},\zeta^\alpha)]\zeta^{\alpha-1}\mbox{\rm
d}\zeta+\mbox{\rm Re}w_0(z).\end{array}\eeq

By Lemma 2.2 and \eqref{3.52}, we obtain, for all
$t\in
\partial\Omega,$  $\lim\limits_{z\in \Omega,\ z\rightarrow
t}\mbox{\rm Re}w(z)=\gamma(t)$. Thus, the proof is completed.

\bre
When $\alpha=1$, Theorem 2.2 is just the one result in \cite{B2009}.
\eere

\section{Schwarz problem for inhomogeneous polyanalytic equations}

 To solve Schwarz problem for the inhomogeneous polyanalytic equation, we define a
poly-Schwarz operator and a Pompeiu operator for $\Omega$ as follows.
\begin{equation}\label{52.2}\begin{array}{ll}&S_n[\gamma_0,\gamma_1,\cdots,\gamma_{n-1}](z)\\[3mm]&=\displaystyle\sum\limits_{k=0}^{n-1}\frac{(-1)^k}{k!}\bigg\{\frac{\alpha}{2\pi
i}\int_{\Gamma_1\bigcup\Gamma_2}\gamma_k(\zeta)\left(\zeta-z+\overline{\zeta-z}\right)^kH_1(z^\alpha,\zeta^\alpha)\displaystyle\frac{\mbox{\rm
d}\zeta}{\zeta}\\[7mm]
&\ \ \ \ +\displaystyle\frac{\alpha}{\pi
i}\int_{[\varpi,w]\cup[r,1]}\gamma_k(\zeta)\left(\zeta-z+\overline{\zeta-z}\right)^kH_2(z^\alpha,\zeta^\alpha)\zeta^{\alpha-1}\mbox{\rm
d}\zeta\bigg\}, \ z\in\Omega,
\end{array}\end{equation}
\begin{equation}\label{55.4}\begin{array}{ll}
T_n[f](z)=\displaystyle\frac{(-1)^n\alpha}{\pi
(n-1)!}&\displaystyle\iint_{\Omega}\left(\zeta-z+\overline{\zeta-z}\right)^{n-1}\bigg[\zeta^{\alpha-1}H_2(z^\alpha,\zeta^\alpha)f(\zeta)\\[3mm]
&\quad\quad-\displaystyle\overline{\zeta^{\alpha-1}}H_2(z^\alpha,\overline{\zeta^\alpha})\overline{f(\zeta)}\bigg]\mbox{\rm
d}\xi\mbox{\rm d}\eta,\ \ z\in\Omega
\end{array}\end{equation}
with $f\in L_p(\Omega;C),\ p>2$.

The operators $T_n$ and $S_n$ have the following differentiability and boundary properties.

\ble When $\gamma_0,\gamma_1,\cdots,\gamma_{n-1}\in
C(\partial\Omega,\mathbb{R})$ and $t\in
\partial\Omega$, then
$$\left\{\begin{array}{ll}\displaystyle\frac{\partial^nS_n[\gamma_0,\gamma_1,\cdots,\gamma_{n-1}](z)}{\partial\overline{z}^n}
=0,&z\in\Omega,\\[3mm]
\left\{\mbox{\rm
Re}\displaystyle\frac{\partial^kS_n[\gamma_0,\gamma_1,\cdots,\gamma_{n-1}]}
{\partial\overline{z}^k}\right\}^+(t)=\gamma_k(t),&k=0,1,2,\cdots,n-1\end{array}\right.$$
with $\displaystyle\frac{\partial^kS_n[\gamma_0,\gamma_1,\cdots,\gamma_{n-1}](z)}{\partial\overline{z}^k}=S_n[\gamma_0,\gamma_1,\cdots,\gamma_{n-1}](z)$ for $z\in \Omega$ and $k=0$.
\eele

\p  Obviously, according to \eqref{52.2}, the first equation in this lemma is true. Furthermore, for $k=0,1,2,\cdots,n-1$,
$$\begin{array}{ll}
&\displaystyle\frac{\partial^kS_n[\gamma_0,\gamma_1,\cdots,\gamma_{n-1}](z)}{\partial\overline{z}^k}\\[3mm]
&=\displaystyle\sum\limits_{l=k}^{n-1}\frac{(-1)^{l-k}}{(l-k)!}\bigg\{\frac{\alpha}{2\pi
i}\int_{\Gamma_1\bigcup \Gamma_2}\gamma_l(\zeta)\left(\zeta-z+\overline{\zeta-z}\right)^{l-k}H_1(z^\alpha,\zeta^\alpha)\displaystyle\frac{\mbox{\rm
d}\zeta}{\zeta}\\[5mm]
&\ \ \ \ +\displaystyle\frac{\alpha}{\pi
i}\int\limits_{[\varpi,\omega]\cup[r,1]}\gamma_l(\zeta)\left(\zeta-z+\overline{\zeta-z}\right)^{l-k}H_2(z^\alpha,\zeta^\alpha)\zeta^{\alpha-1}\mbox{\rm
d}\zeta\bigg\}\\[7mm]
&=\displaystyle\sum\limits_{l=k}^{n-1}\sum\limits_{j=0}^{l-k}\frac{(-1)^{l-k}}{j!(l-k-j)!}(-z-\overline{z})^{l-k-j}
\bigg\{\displaystyle\frac{\alpha}{2\pi
i}\int_{\Gamma_1\bigcup \Gamma_2}\gamma_l(\zeta)\left(\zeta+\overline{\zeta}\right)^{j}H_1(z^\alpha,\zeta^\alpha)\displaystyle\frac{\mbox{\rm
d}\zeta}{\zeta}\\[5mm]
&\ \ \ \  +\displaystyle\frac{\alpha}{\pi
i}\int_{[\varpi,\omega]\cup[r,1]}\gamma_l(\zeta)\left(\zeta+\overline{\zeta}\right)^{j}H_2(z^\alpha,\zeta^\alpha)\zeta^{\alpha-1}\mbox{\rm
d}\zeta\bigg\}.\end{array}$$

 By simple calculation, for $\zeta\in \Gamma_1\bigcup\Gamma_2$,
$$
\displaystyle\frac{1}{2\zeta}\left[H_1(z^\alpha,\zeta^\alpha)+H_1(\overline{z^\alpha},\overline{\zeta^\alpha})\right]=
\left[H_2(z^\alpha,\zeta^\alpha)-H_2(\overline{z^\alpha},\zeta^\alpha)\right]\zeta^{\alpha-1},
$$
which implies
$$\begin{array}{ll}\mbox{\rm
Re}\left\{\displaystyle\frac{\partial^kS_n[\gamma_0,\gamma_1,\cdots,\gamma_{n-1}](z)}{\partial\overline{z}^k}\right\}
\\[5mm]=\displaystyle\sum\limits_{l=k}^{n-1}\sum_{j=0}^{l-k}\frac{(-1)^{l-k}(-z-\overline{z})^{l-k-j}}{j!(l-k-j)!}\frac{\alpha}{2\pi
i}\int_{\partial\Omega}\gamma_l(\zeta)\left(\zeta+\overline{\zeta}\right)^{j}G(z,\zeta)\zeta^{\alpha-1}{\mbox{\rm
d}}\zeta,\end{array}$$
where
\begin{equation}\label{55.1}
G(z,\zeta)=H_2(z^\alpha,\zeta^\alpha)-H_2(\overline{z}^\alpha,\zeta^\alpha),\ \ z,\,\zeta\in \Omega.\end{equation}

By Lemma 2.2, $$ \lim\limits_{z\rightarrow t,t\in
\partial\Omega}\frac{\alpha}{2\pi
i}\int_{\partial\Omega}\gamma_l(\zeta)\left(\zeta+\overline{\zeta}\right)^{j}G(z,\zeta)\zeta^{\alpha-1}{\mbox{\rm
d}}\zeta =\gamma_l(t)\left(t+\overline{t}\right)^{j},$$
hence,
$$
\left\{\mbox{\rm
Re}\displaystyle\frac{\partial^kS_n[\gamma_0,\gamma_1,\cdots,\gamma_{n-1}]}{\partial\overline{z}^k}\right\}^+(t)
=\displaystyle\sum\limits_{l=k}^{n-1}\sum\limits_{j=0}^{l-k}\frac{(-1)^{j}(t+\overline{t})^{l-k}\gamma_l(\zeta)}{j!(l-k-j)!}
=\gamma_k(t).$$
Therefore, the proof is finished.

\ble For $f\in L_p(\Omega;\mathbb{C}),\ p>2$, $z\in\Omega$ and $t\in\partial\Omega$,
$$\left\{\begin{array}{ll} \displaystyle\frac{\partial
T_k[f](z)}{\partial\overline{z}}=T_{k-1}[f](z), &k\geq2,\,\,\, k\in\mathbb{N},\\[3mm]
\displaystyle\frac{\partial^nT_n[f](z)}{\partial\overline{z}^n}=f(z),&n\geq1,\,\,\,
n\in\mathbb{N},\\[3mm]
 \left\{\mbox{\rm Re}\partial^{k}_{\overline{z}}T_n[f]\right\}^+(t)=0,& k=0,1,\cdots,n-1.\end{array}\right.$$
\eele

\p When $k\geq2,$ we know
\begin{equation}\label{2.6}
T_k[f](z)=\displaystyle\frac{(-1)^{k-1}}{(k-1)!}\sum\limits_{j=0}^{k-1}\left(\begin{array}{cc}k-1\\j\end{array}\right)
(-z-\overline{z})^{k-j-1}T_1[(\zeta+\overline{\zeta})^jf](z)
\end{equation}
with
\begin{equation}
T_1[f](z)=-\displaystyle\frac{\alpha}{\pi}\displaystyle\iint_{\Omega}\left[\zeta^{\alpha-1}H_2(z^\alpha,\zeta^\alpha)f(\zeta)-\displaystyle\overline{\zeta^{\alpha-1}}H_2(z^\alpha,\overline{\zeta^\alpha})\overline{f(\zeta)}\right]\mbox{\rm
d}\xi\mbox{\rm d}\eta.\end{equation}
From the proof in Theorem 2.2,
$\displaystyle\frac{\partial
T_1[f](z)}{\partial\overline{z}}=f(z)$. Then
$$\begin{array}{ll}
\displaystyle\frac{\partial
T_k[f](z)}{\partial\overline{z}}&=\displaystyle\frac{(-1)^{k-2}}{(k-2)!}\sum\limits_{j=0}^{k-2}
\left(\begin{array}{cc}k-2\\j\end{array}\right)(-z-\overline{z})^{k-j-2}T_1[(\zeta+\overline{\zeta})^jf](z)\\[5mm]
&=T_{k-1}[f](z),
\end{array}$$
which implies
$\displaystyle\frac{\partial^nT_n[f](z)}{\partial\overline{z}^n}=\displaystyle\frac{\partial
T_1[f](z)}{\partial\overline {z}}=f(z)$ for $n\geq1$. Furthermore, by
\eqref{2.6}, we obtain for $n\geq1$,
\beq\label{9990000}\mbox{\rm
Re}\{T_n[f](z)\}=\displaystyle\frac{(-1)^{n-1}}{(n-1)!}\sum\limits_{k=0}^{n-1}
\left(\begin{array}{cc}n-1\\k\end{array}\right)(-z-\overline{z})^{n-k-1}\mbox{\rm
Re}\left\{T_1[(\zeta+\overline{\zeta})^kf](z)\right\},\eeq
where
$$
\mbox{\rm
Re}\left\{T_1[(\zeta+\overline{\zeta})^kf](z)\right\}=-\displaystyle\frac{\alpha}{2\pi}\iint_{\Omega}(\zeta+\overline{\zeta})^k\left[f(\zeta)G^*(z,\zeta)+\overline{f(\zeta)}\
\overline{G^*(z,\zeta)}\right]\mbox{\rm d}\xi \mbox{\rm d}\eta$$
with $G^*(z,\zeta)$ given in the proof of Theorem 2.2, and for $(z,\zeta)\in
\partial\Omega\times\Omega$, $G^*(z,\zeta)=0$. Thus, $\mbox{\rm
Re}\left\{T_1[(\zeta+\overline{\zeta})^kf]\right\}^+(t)\}=0$ for $t\in\partial\Omega$, which implies $\{\mbox{\rm
Re}T_n[f]\}^+(t)=0$ for $n\geq1$. In addition, from the first equation of this lemma, $\displaystyle\frac{\partial^k
T_n[f](z)}{\partial\overline{z}^k}=T_{n-k}[f](z)$ for $k=0,1,\cdots,n-1$, then  $\left\{\mbox{\rm Re}\partial^{k}_{\overline{z}}T_n[f]\right\}^+(t)=0$ for $t\in\partial\Omega$.
Therefore, the proof is completed.

By Lemmas 3.1-3.2, we can solve the following Schwarz problem for polyanalytic equation.
\bth The Schwarz problem for polyanalytic equation in
$\Omega$ \begin{equation}\label{2.7} \left\{\begin{array}{ll}
\partial^n_{\overline{z}}w(z)=f(z), &
z\in\Omega,\ \ f\in L_p(\Omega,\mathbb{C}),\ p>2,\\[2mm]
\left\{\mbox{\rm
Re}(\partial^k_{\overline{z}}w)\right\}^+(t)=\gamma_k(t),&t\in\partial\Omega,\
\gamma_k\in C(\partial\Omega,\ \mathbb{R}),\ k=0,1,\cdots,n-1,
\end{array}\right.
\end{equation}
is solvable by
$$
w(z)=S_n[\gamma_0,\gamma_1,\cdots,\gamma_{n-1}](z)+T_n[f](z)
+i\sum\limits_{k=0}^{n-1}(z+\overline{z})^kc_k,
$$
where $c_k\in \mathbb{R}$, $S_n,\ T_n$ are defined by \eqref{52.2}
and \eqref{55.4}, respectively.
\eeth

\p  From Lemmas 3.1 and 3.2, it is easy to know that $S_n[\gamma_0,\gamma_1,\cdots,\gamma_{n-1}](z)+T_n[f](z)
$ satisfies the condition \eqref{2.7}. So we write
$$w(z)=S_n[\gamma_0,\gamma_1,\cdots,\gamma_{n-1}](z)+T_n[f](z)+U(z),$$
which means
\beq\label{000111}
\left\{\begin{array}{ll} \partial^n_{\overline{z}}U(z)=0, &
z\in\Omega;\\[2mm]
\left\{\mbox{\rm
Re}(\partial^k_{\overline{z}}U)\right\}^+(t)=0,&t\in\partial\Omega,\
k=0,1,\cdots,n-1.
\end{array}\right.
\eeq
According to the result in \cite{bb}, the polyanalytic function $U(z)$ of order $n$ has the expression as follows
\begin{equation}\label{2.8}
U(z)=\sum\limits_{k=0}^{n-1}\displaystyle(z+\overline{z})^kf_k(z),\
\ \ z\in \Omega
\end{equation}
with $f_k$ being analytic functions. Taking \eqref{2.8} into the
second equation of \eqref{000111}, $\{\mbox{\rm Re}f_k\}^+(t)=0,\ \ t\in
\partial\Omega$. Thus, by the Theorem 2.1, $f_k(z)=ic_k$ where $c_k$ are real numbers. The proof is completed.

\section{Schwarz problem for a generalized polyanalytic equation}

In this section, we consider Schwarz problem for a generalized polyanalytic equation for the domain $R^+=\{\,0<r<|z|<1, \ \mbox{\rm Im}z>0\}$. We firstly discuss some properties in more detail for the operator $T_n$, which will be used in the sequel. Let
\begin{equation}\label{6.1}\begin{array}{ll}
\widetilde{T}_k[f](z)=\displaystyle\frac{(-1)^k}{\pi
(k-1)!}&\displaystyle\iint_{R^+}\left(\zeta-z+\overline{\zeta-z}\right)^{k-1}\bigg[H(z,\zeta)f(\zeta)\\[3mm]
&\quad\quad-\displaystyle H(z,\overline{\zeta})\overline{f(\zeta)}\bigg]\mbox{\rm
d}\xi\mbox{\rm d}\eta,\ \ z\in R^+,
\end{array}\end{equation}
which is just \eqref{55.4} for $\alpha=1$. Here
\begin{equation}\label{6.12*}\begin{array}{ll}
H(z,\zeta)&=\displaystyle\frac{1}{\zeta-z}-\displaystyle\frac{z}{1-z\zeta}\\[3mm]&+\displaystyle\sum\limits_{n=1}^\infty
r^{2n}\left[\displaystyle\frac{1}{r^{2n}\zeta-z}-\displaystyle\frac{z}{\zeta(r^{2n}z-\zeta)}-\displaystyle\frac{1}{\zeta(r^{2 n}-z\zeta)}
+\displaystyle\frac{z}{r^{2 n}z\zeta-1}\right].\end{array}\end{equation}
Then by Lebiniz rule, we obtain for $0\leq l\leq k-1$
\beq\label{6.3}\begin{array}{ll}
\displaystyle\frac{\partial^l
\widetilde{T}_k[f](z)}{\partial z^l}
=\displaystyle\frac{1}{\pi}\displaystyle\iint_{R^+}\displaystyle\sum\limits_{j=0}^l&\left(\begin{array}{cc}1\\j\end{array}\right)\displaystyle\frac{(-1)^{k-j}(l-j)!}{(k-j-1)!}
\left(\zeta-z+\overline{\zeta-z}\right)^{k-j-1}\\[5mm]&\left[f(\zeta)g(z,\zeta)-\overline{f(\zeta)}g(z,\overline{\zeta})\right]
\mbox{\rm
d}\xi\mbox{\rm d}\eta,
\end{array}
\eeq
where
\beq\label{6.4}
\begin{array}{ll}
g(z,\zeta)=\displaystyle\frac{1}{(\zeta-z)^{l-j+1}}-\displaystyle\frac{\zeta^{l-j-1}}{(1-z\zeta)^{l-j+1}}+\sum\limits_{n=1}^\infty
\bigg[\displaystyle\frac{r^{2n}}{(r^{2n}\zeta-z)^{l-j+1}}\\[5mm]\quad-\displaystyle\frac{(-1)^{l-j}r^{2n(l-j)}}{(r^{2n}z-\zeta)^{l-j+1}}-
\displaystyle\frac{r^{2n}\zeta^{l-j-1}}{(r^{2n}-z\zeta)^{l-j+1}}
+\displaystyle\frac{(-1)^{l-j}r^{2n(l-j)}\zeta^{l-j-1}}{(r^{2n}z\zeta-1)^{l-j+1}}\bigg].
\end{array}\eeq

\ble
Let $f\in L_p(R^+,\mathbb{C})$, $p>2$ and $k\in \mathbb{N}$, then
$$
\left|\displaystyle\frac{\partial^l
\widetilde{T}_k[f](z)}{\partial z^l}\right|\leq C(l,p)\parallel f\parallel_{L_p(R^+)},\ \ z\in \Omega,\ \ l=0,1,\cdots,k-1,
$$
where $C(l,p)$ is a constant depending on $l$ and $p$.
\eele

\p  By simple calculation, for $z$, $\zeta\in R^+$, there exists some $M>0$ such that $|r^{2n}\zeta-z|>M,$ $|r^{2n}z-\zeta|>M$, $|r^{2n}z\zeta-1|>M$ for $n\geq1$, and
$|r^{2n}-z\zeta|>M$ for $n\geq2$. Then $g$ in \eqref{6.4} can be rewritten
\beq\label{2021999}
g(z,\zeta)=\displaystyle\frac{1}{(\zeta-z)^{l-j+1}}-\displaystyle\frac{\zeta^{l-j-1}}{(1-z\zeta)^{l-j+1}}
-\displaystyle\frac{r^{2}\zeta^{l-j-1}}{(r^{2}-z\zeta)^{l-j+1}}+g^*(z,\zeta)
\eeq
with
\beq\label{20212021}\begin{array}{ll}
g^*(z,\zeta)=\displaystyle\sum\limits_{n=1}^\infty
\bigg[\displaystyle\frac{r^{2n}}{(r^{2n}\zeta-z)^{l-j+1}}&-\displaystyle\frac{(-1)^{l-j}r^{2n(l-j)}}{(r^{2n}z-\zeta)^{l-j+1}}-
\displaystyle\frac{r^{2(n+1)}\zeta^{l-j-1}}{(r^{2n+2}-z\zeta)^{l-j+1}}\\[5mm]
&+\displaystyle\frac{(-1)^{l-j}r^{2n(l-j)}\zeta^{l-j-1}}{(r^{2n}z\zeta-1)^{l-j+1}}\bigg],\end{array}
\eeq
which is bounded. Similarly, we have $g^*(z,\overline{\zeta})$ is also bounded for $z,\zeta\in R^+$. In addition, from \eqref{6.3}, we rewrite
\beq\label{6.000}
\displaystyle\frac{\partial^l
\widetilde{T}_k[f](z)}{\partial z^l}
=\displaystyle\sum\limits_{j=0}^l\left(\begin{array}{cc}1\\j\end{array}\right)\displaystyle\frac{(-1)^{k-j}(l-j)!}{(k-j-1)!}\left[A_{l,j}(z)
+B_{l,j}(z)\right],
\eeq
where
\beq\label{6.001}\begin{array}{ll}
A_{l,j}(z)\\[3mm]
=\displaystyle\frac{1}{\pi}\displaystyle\iint_{R^+}\left(\zeta-z+\overline{\zeta-z}\right)^{k-j-1}\Bigg\{
f(\zeta)\left[\displaystyle\frac{1}{(\zeta-z)^{l-j+1}}-\displaystyle\frac{\zeta^{l-j-1}}{(1-z\zeta)^{l-j+1}}
-\displaystyle\frac{r^{2}\zeta^{l-j-1}}{(r^{2}-z\zeta)^{l-j+1}}\right]\\[5mm]\quad-\overline{f(\zeta)}\left[\displaystyle\frac{1}{(\overline{\zeta}-z)^{l-j+1}}
-\displaystyle\frac{\overline{\zeta}^{l-j-1}}{(1-z\overline{\zeta})^{l-j+1}}
-\displaystyle\frac{r^{2}\overline{\zeta}^{l-j-1}}{(r^{2}-z\overline{\zeta})^{l-j+1}}\right]\Bigg\}\mbox{\rm
d}\xi\mbox{\rm d}\eta,\end{array}
\eeq
and
\beq\label{6.002}
B_{l,j}(z)=\displaystyle\frac{1}{\pi}\displaystyle\iint_{R^+}\left(\zeta-z+\overline{\zeta-z}\right)^{k-j-1}\left[
f(\zeta)g^*(z,\zeta)-\overline{f(\zeta)}g^*(z,\overline{\zeta})\right]\mbox{\rm
d}\xi\mbox{\rm d}\eta.
\eeq

Since for $z$, $\zeta$ in $R^+$, we have the fact that
\beq\label{888}|\zeta-z|<|\overline{\zeta}-z|, \,\,\,|\zeta-z|<|1-z\overline{\zeta}|<|1-z\zeta|,\eeq
and
 \beq\label{999}r|\zeta-z|<|r^2-z\overline{\zeta}|<|r^2-z\zeta|,\eeq thus, making use of H$\ddot{o}$lder inequality, for $\displaystyle\frac{1}{p}+\displaystyle\frac{1}{q}=1$,
$$
|A_{l,j}(z)|\leq C(l,j)\left(\displaystyle\iint_{R^+}\left|\zeta-z\right|^{q(k-l-2)}\mbox{\rm
d}\xi\mbox{\rm d}\eta\right)^{1/q}\parallel f\parallel_{L_p(R^+)}, \quad  j=0,1,\cdots,l.
$$
Therefore,
$$|A_{l,j}(z)|\leq C(l,p,j)\parallel f \parallel_{L_p(R^+)}$$
is true since $k-l-2\geq-1$ and $p>2$, where $C(l,p,j)$ is a constant depending on $l,p,j$.
Also,
$$
|B_{l,j}(z)|\leq C(l,j)\left(\displaystyle\iint_{R^+}\left|\zeta-z\right|^{q(k-j-1)}\mbox{\rm
d}\xi\mbox{\rm d}\eta\right)^{1/q}\parallel f\parallel_{L_p(R^+)}, \quad  j=0,1,\cdots,l,
$$
 and $k-j-1\geq 0$, which implies $|B_{l,j}(z)|\leq C(l,p,j)\parallel f \parallel_{L_p(R^+)}$. Therefore, the proof is finished.

\ble
Suppose $f\in L_p(R^+)$, $p>2$, then for $z_1,z_2\in R^+$ and $k\in \mathbb{N}$,
$$
\left|
\displaystyle\frac{\partial^l
\widetilde{T}_k[f](z_1)}{\partial z^l}-\displaystyle\frac{\partial^l
\widetilde{T}_k[f](z_2)}{\partial z^l}
\right|\leq C(l,p)\parallel f\parallel_{L_p(R^+)}\left\{\begin{array}{ll} |z_1-z_2|,&0\leq l\leq k-2,\\[3mm]
|z_1-z_2|^{(p-2)/p},& l= k-1,\end{array}\right.
$$
where $C(l,p)$ is a constant depending on $l$ and $p$.
\eele

\p When $l=0,1,\cdots,k-2$, from Lemma 4.1, we know
$$
\left|\displaystyle\frac{\partial^{l+1}
\widetilde{T}_k[f](z)}{\partial z^{l+1}}\right|\leq C(l,p)\parallel f\parallel_{L_p(R^+)},
$$
then by the mean value theorem,
\beq\label{20210}
\left|
\displaystyle\frac{\partial^l
\widetilde{T}_k[f](z_1)}{\partial z^l}-\displaystyle\frac{\partial^l
\widetilde{T}_k[f](z_2)}{\partial z^l}
\right|\leq C(l,p)\parallel f\parallel_{L_p(R^+)}|z_1-z_2|.
\eeq
For $l=k-1$, from \eqref{6.000},
\beq\label{20211}\begin{array}{ll}
\left|
\displaystyle\frac{\partial^{k-1}
\widetilde{T}_k[f](z_1)}{\partial z^{k-1}}-\displaystyle\frac{\partial^{k-1}
\widetilde{T}_k[f](z_2)}{\partial z^{k-1}}
\right|\\[5mm]\leq \displaystyle\sum\limits_{j=0}^{k-1}\left(\begin{array}{cc}k-1\\j\end{array}\right)\left[\left|A_{k-1,j}(z_1)-A_{k-1,j}(z_2)\right|+
\left|B_{k-1,j}(z_1)-B_{k-1,j}(z_2)\right|
\right]
\end{array}\eeq
with $A_{k-1,j},$ $B_{k-1,j}$ given by \eqref{6.001} and \eqref{6.002} for $l=k-1$, respectively.

Invoking the identity
\beq\label{000}
b^{i+1}c^i-a^{i+1}d^i=(b-a)\displaystyle\sum\limits_{j=0}^ib^{i-j}a^jc^{i-j}d^j+(c-d)\displaystyle\sum\limits_{j=0}^{i-1}b^{i-j}a^{j+1}c^{i-1-j}d^j,\eeq
and $|\zeta-z|\leq |1-z\zeta|$ for $z,\zeta\in R^+$, we obtain for $j=0,1,\cdots,k-1,$
$$\begin{array}{ll}
\left|(1-z_2\zeta)^{k-j}\left(\zeta-z_1+\overline{\zeta-z_1}\right)^{k-j-1}-(1-z_1\zeta)^{k-j}\left(\zeta-z_2+\overline{\zeta-z_2}\right)^{k-j-1}\right|\\[5mm]
\leq C(k)|z_1-z_2||1-z_1\zeta|^{k-j-1}|1-z_2\zeta|^{k-j-1},\end{array}$$
where $C(k)$ is a constant depending on $l$, which gives
$$\begin{array}{ll}
\left|\displaystyle\frac{\left(\zeta-z_1+\overline{\zeta-z_1}\right)^{k-j-1}}{(1-z_1\zeta)^{k-j}}
-\displaystyle\frac{\left(\zeta-z_2+\overline{\zeta-z_2}\right)^{k-j-1}}{(1-z_2\zeta)^{k-j}}\right|\\[5mm]
=\displaystyle\frac{\left|(1-z_2\zeta)^{k-j}\left(\zeta-z_1+\overline{\zeta-z_1}\right)^{k-j-1}-(1-z_1\zeta)^{k-j}\left(\zeta-z_2+\overline{\zeta-z_2}\right)^{k-j-1}\right|}{|1-z_1\zeta|^{k-j}|1-z_2\zeta|^{k-j}}\\[5mm]
\leq \displaystyle\frac{C(k)|z_1-z_2|}{\left|1-z_2\zeta\right|\left|1-z_1\zeta\right|}
\leq\displaystyle\frac{C(k)|z_1-z_2|}{\left|\zeta-z_1\right|\left|\zeta -z_2\right|},
\end{array}
$$

Similarly, for $z,\zeta\in R^+$, in terms of \eqref{888},\eqref{999} and \eqref{000}, we obtain
$$
\left|\displaystyle\frac{\left(\zeta-z_1+\overline{\zeta-z_1}\right)^{k-j-1}}{(\zeta-z_1)^{k-j}}
-\displaystyle\frac{\left(\zeta-z_2+\overline{\zeta-z_2}\right)^{k-j-1}}{(\zeta-z_2)^{k-j}}\right|
\leq\displaystyle\frac{C(k)|z_1-z_2|}{\left|\zeta-z_1\right|\left|\zeta -z_2\right|},$$
$$
\left|\displaystyle\frac{\left(\zeta-z_1+\overline{\zeta-z_1}\right)^{k-j-1}}{(r^{2}-z_1\zeta)^{k-j}}
-\displaystyle\frac{\left(\zeta-z_2+\overline{\zeta-z_2}\right)^{k-j-1}}{(r^{2}-z_2\zeta)^{k-j}}\right|
\leq\displaystyle\frac{C(k)|z_1-z_2|}{\left|\zeta-z_1\right|\left|\zeta -z_2\right|},$$
$$
\left|\displaystyle\frac{\left(\zeta-z_1+\overline{\zeta-z_1}\right)^{k-j-1}}{(\overline{\zeta}-z_1)^{k-j}}
-\displaystyle\frac{\left(\zeta-z_2+\overline{\zeta-z_2}\right)^{k-j-1}}{(\overline{\zeta}-z_2)^{k-j}}\right|
\leq\displaystyle\frac{C(k)|z_1-z_2|}{\left|\zeta-z_1\right|\left|\zeta -z_2\right|},$$
$$
\left|\displaystyle\frac{\left(\zeta-z_1+\overline{\zeta-z_1}\right)^{k-j-1}}{(1-z_1\overline{\zeta})^{k-j}}
-\displaystyle\frac{\left(\zeta-z_2+\overline{\zeta-z_2}\right)^{k-j-1}}{(1-z_2\overline{\zeta})^{k-j}}\right|
\leq\displaystyle\frac{C(k)|z_1-z_2|}{\left|\zeta-z_1\right|\left|\zeta -z_2\right|},$$
$$
\left|\displaystyle\frac{\left(\zeta-z_1+\overline{\zeta-z_1}\right)^{k-j-1}}{(r^2-z_1\overline{\zeta})^{k-j}}
-\displaystyle\frac{\left(\zeta-z_2+\overline{\zeta-z_2}\right)^{k-j-1}}{(r^2-z_2\overline{\zeta})^{k-j}}\right|
\leq\displaystyle\frac{C(k)|z_1-z_2|}{\left|\zeta-z_1\right|\left|\zeta -z_2\right|},$$
thus,
$$
\left|
A_{k-1,j}(z_1)-A_{k-1,j}(z_2)\right|
\leq C(k)|z_1-z_2|\displaystyle\iint_{R^+}\frac{|f(\zeta)|}{\left|\zeta-z_1\right|\left|\zeta -z_2\right|}\mbox{\rm
d}\xi\mbox{\rm d}\eta.
$$
Then,
\beq\label{202000}
\left|
A_{k-1,j}(z_1)-A_{k-1,j}(z_2)\right|
\leq C(k,p)|z_1-z_2|^{(p-2)/p}\parallel f\parallel_{L_p(R^+)}.
\eeq

By \eqref{000}, for $z,\zeta\in R^+$
$$\begin{array}{ll}
\left|(r^{2n}\zeta-z_2)^{k-j}\left(\zeta-z_1+\overline{\zeta-z_1}\right)^{k-j-1}-(r^{2n}\zeta-z_1)^{k-j}\left(\zeta-z_2+\overline{\zeta-z_2}\right)^{k-j-1}\right|\\[5mm]
\leq C(k)|z_1-z_2|.\end{array}$$
Also, for $z$, $\zeta\in R^+$, there exists some $M>0$ such that $|r^{2n}\zeta-z|>M,$ therefore,
$$\left|\displaystyle\frac{\left(\zeta-z_1+\overline{\zeta-z_1}\right)^{k-j-1}}{(r^{2n}\zeta-z_1)^{k-j}}
-\displaystyle\frac{\left(\zeta-z_2+\overline{\zeta-z_2}\right)^{k-j-1}}{(r^{2n}\zeta-z_2)^{k-j}}\right|
\leq C(k)|z_1-z_2|,$$
which implies
$$\begin{array}{ll}
\left|\displaystyle\sum\limits_{n=1}^{\infty}r^{2n}\left[\displaystyle\frac{\left(\zeta-z_1+\overline{\zeta-z_1}\right)^{k-j-1}}{(r^{2n}\zeta-z_1)^{k-j}}
-\displaystyle\frac{\left(\zeta-z_2+\overline{\zeta-z_2}\right)^{k-j-1}}{(r^{2n}\zeta-z_2)^{k-j}}\right]\right|\\[5mm]
\leq C(k)|z_1-z_2|\displaystyle\sum\limits_{n=1}^{\infty}r^{2n}
\leq C(k)|z_1-z_2|.
\end{array}
$$
Similarly,
$$\begin{array}{ll}
\left|\displaystyle\sum\limits_{n=1}^{\infty}(-r^{2n})^{k-j-1}\left[\displaystyle\frac{\left(\zeta-z_1+\overline{\zeta-z_1}\right)^{k-j-1}}{(r^{2n}z_1-\zeta)^{k-j}}
-\displaystyle\frac{\left(\zeta-z_2+\overline{\zeta-z_2}\right)^{k-j-1}}{(r^{2n}z_2-\zeta)^{k-j}}\right]\right|
\leq C(k)|z_1-z_2|,
\end{array}
$$
$$\begin{array}{ll}
\left|\displaystyle\sum\limits_{n=1}^{\infty}r^{2(n+1)}\left[\displaystyle\frac{\left(\zeta-z_1+\overline{\zeta-z_1}\right)^{k-j-1}}{(r^{2n+2}-z_1\zeta)^{k-j}}
-\displaystyle\frac{\left(\zeta-z_2+\overline{\zeta-z_2}\right)^{k-j-1}}{(r^{2n+2}-z_2\zeta)^{k-j}}\right]\right|
\leq C(k)|z_1-z_2|,
\end{array}
$$
$$\begin{array}{ll}
\left|\displaystyle\sum\limits_{n=1}^{\infty}(-r^{2n})^{k-j-1}\left[\displaystyle\frac{\left(\zeta-z_1+\overline{\zeta-z_1}\right)^{k-j-1}}{(r^{2n}z_1\zeta-1)^{k-j}}
-\displaystyle\frac{\left(\zeta-z_2+\overline{\zeta-z_2}\right)^{k-j-1}}{(r^{2n}z_2\zeta-1)^{k-j}}\right]\right|
\leq C(k)|z_1-z_2|,
\end{array}
$$
and
$$\begin{array}{ll}
\left|\displaystyle\sum\limits_{n=1}^{\infty}r^{2n}\left[\displaystyle\frac{\left(\zeta-z_1+\overline{\zeta-z_1}\right)^{k-j-1}}{(r^{2n}\overline{\zeta}-z_1)^{k-j}}
-\displaystyle\frac{\left(\zeta-z_2+\overline{\zeta-z_2}\right)^{k-j-1}}{(r^{2n}\overline{\zeta}-z_2)^{k-j}}\right]\right|
\leq C(k)|z_1-z_2|,
\end{array}
$$
$$\begin{array}{ll}
\left|\displaystyle\sum\limits_{n=1}^{\infty}(-r^{2n})^{k-j-1}\left[\displaystyle\frac{\left(\zeta-z_1+\overline{\zeta-z_1}\right)^{k-j-1}}{(r^{2n}z_1-\overline{\zeta})^{k-j}}
-\displaystyle\frac{\left(\zeta-z_2+\overline{\zeta-z_2}\right)^{k-j-1}}{(r^{2n}z_2-\overline{\zeta})^{k-j}}\right]\right|
\leq C(k)|z_1-z_2|,
\end{array}
$$
$$\begin{array}{ll}
\left|\displaystyle\sum\limits_{n=1}^{\infty}r^{2(n+1)}\left[\displaystyle\frac{\left(\zeta-z_1+\overline{\zeta-z_1}\right)^{k-j-1}}{(r^{2n+2}-z_1\overline{\zeta})^{k-j}}
-\displaystyle\frac{\left(\zeta-z_2+\overline{\zeta-z_2}\right)^{k-j-1}}{(r^{2n+2}-z_2\overline{\zeta})^{k-j}}\right]\right|
\leq C(k)|z_1-z_2|,
\end{array}
$$
$$\begin{array}{ll}
\left|\displaystyle\sum\limits_{n=1}^{\infty}(-r^{2n})^{k-j-1}\left[\displaystyle\frac{\left(\zeta-z_1+\overline{\zeta-z_1}\right)^{k-j-1}}{(r^{2n}z_1\overline{\zeta}-1)^{k-j}}
-\displaystyle\frac{\left(\zeta-z_2+\overline{\zeta-z_2}\right)^{k-j-1}}{(r^{2n}z_2\overline{\zeta}-1)^{k-j}}\right]\right|
\leq C(k)|z_1-z_2|.
\end{array}
$$
Thus, from the expression of $g^*(z,\zeta)$ in \eqref{20212021} for $l=k-1$, we obtain
$$
\left|
\left(\zeta-z_1+\overline{\zeta-z_1}\right)^{k-j-1}g^*(z_1,\zeta)-\left(\zeta-z_2+\overline{\zeta-z_2}\right)^{k-j-1}g^*(z_2,\zeta)\right|
\leq C(k)|z_1-z_2|,
$$
and
$$
\left|
\left(\zeta-z_1+\overline{\zeta-z_1}\right)^{k-j-1}g^*(z_1,\overline{\zeta})-\left(\zeta-z_2+\overline{\zeta-z_2}\right)^{k-j-1}g^*(z_2,\overline{\zeta})\right|
\leq C(k)|z_1-z_2|,
$$
which means
\beq\label{202011}
\left|
B_{k-1,j}(z_1)-B_{k-1,j}(z_2)\right|
\leq C(k,p)|z_1-z_2|\parallel f\parallel_{L_p(R^+)}.
\eeq
Therefore, from the \eqref{20211},\eqref{202000} and \eqref{202011},
$$
\left|
\displaystyle\frac{\partial^{k-1}
\widetilde{T}_k[f](z_1)}{\partial z^{k-1}}-\displaystyle\frac{\partial^{k-1}
\widetilde{T}_k[f](z_2)}{\partial z^{k-1}}
\right|\leq C(k,p)\parallel f\parallel_{L_P(R^+)}|z_1-z_2|^{(p-2)/p}.$$
That is, the proof is completed.

For the sake of investigating the boundedness of $\displaystyle\frac{\partial^k
\widetilde{T}_k f}{\partial z^k}$ , the following operators are introduced.
\beq\label{202109}
\Xi_1f(z)=\displaystyle\frac{(-1)^k k}{\pi}\displaystyle\iint_{R^+}\left[\displaystyle\frac{(\overline{\zeta}-\overline{z})^{k-1}}{(\zeta-z)^{k+1}}f(\zeta)-
\displaystyle\frac{(\zeta-\overline{z})^{k-1}}{(\overline{\zeta}-z)^{k+1}}\overline{f(\zeta)}\right]\mbox{\rm
d}\xi\mbox{\rm d}\eta,
\eeq
\beq\label{202110}
\Xi_2f(z)=\displaystyle\frac{(-1)^{k+1} k}{\pi}\displaystyle\iint_{R^+}\left\{\displaystyle\frac{\left[|\zeta|^2-1+\zeta(\zeta-\overline{z})\right]^{k-1}f(\zeta)}{(1-z\zeta)^{k+1}}-
\displaystyle\frac{\left[|\zeta|^2-1+\overline{\zeta}(\overline{\zeta}-\overline{z})\right]^{k-1}\overline{f(\zeta)}}{(1-z\overline{\zeta})^{k+1}}\right\}\mbox{\rm
d}\xi\mbox{\rm d}\eta,
\eeq
\beq\label{202111}
\Xi_3f(z)=\displaystyle\frac{(-1)^{k+1} kr^2}{\pi}\displaystyle\iint_{R^+}\left\{\displaystyle\frac{\left[|\zeta|^2-r^2+\zeta(\zeta-\overline{z})\right]^{k-1}f(\zeta)}{(r^2-z\zeta)^{k+1}}-
\displaystyle\frac{\left[|\zeta|^2-r^2+\overline{\zeta}(\overline{\zeta}-\overline{z})\right]^{k-1}\overline{f(\zeta)}}{(r^2-z\overline{\zeta})^{k+1}}\right\}\mbox{\rm
d}\xi\mbox{\rm d}\eta.
\eeq

\ble For $f\in L_p(R^+)$ and $p>2$,
$$\parallel \Xi_1f\parallel_{L_p(R^+)}\leq C(p)\parallel f\parallel_{L_p(R^+)},$$
where $\Xi_1f$ is defined by \eqref{202109} and $C(p)$ is a constant depending on $p$.
\eele

\p
Let $\Xi_1(z)=Tf(z)-\widetilde{T}f(z),$ where
$$Tf(z)=\displaystyle\frac{(-1)^k k}{\pi}\displaystyle\iint_{R^+}\displaystyle\left(\frac{\overline{\zeta-z}}{\zeta-z}\right)^{k-1}\frac{f(\zeta)}{(\zeta-z)^2}\mbox{\rm
d}\xi\mbox{\rm d}\eta$$
and
$$\widetilde{T}f(z)=\displaystyle\frac{(-1)^k k}{\pi}\displaystyle\iint_{R^+}\Omega_k(\overline{\zeta}-z)\overline{f(\zeta)}\mbox{\rm
d}\xi\mbox{\rm d}\eta.$$
Here
$$
\Omega_k(z)=\left(\frac{\overline{z}}{z}\right)^{k-1}\frac{1}{z^2}:=\frac{\Lambda(z)}{|z|^2}
$$
with $\Lambda(z)=\left(\frac{\overline{z}}{z}\right)^k.$

From the result in \cite{hile1997}, the singular integral operator $T$ is bounded on $L_p(R^+)$ for $p>2$. Furthermore, similarly in Lemma 2.6 of \cite{ay2018}, $\Lambda(z)$ is homogeneous of degree zero and
$\displaystyle\int_{|z|=1}\Lambda(z)d\sigma(z)=0,$
where $d\sigma(z)$ is the arc length differential on $|z|=1$. Thus, the operator $\widetilde{T}$ is bounded on $L_p(R^+)$ for $p>2$, which means that $\Xi_1f$ is bounded on $L_p(R^+)$ for $p>2$.

\ble For $f\in L_p(R^+)$ and $p>2$,
$$\parallel \Xi_2f\parallel_{L_p(R^+)}\leq C(p)\parallel f\parallel_{L_p(R^+)},$$
\eele
where $\Xi_2f$ is defined by \eqref{202110} and $C(p)$ is a constant depending on $p$.

\p
According to \eqref{202110}, putting  $\Xi_2(z)=P_1f(z)+P_2f(z)$ with
$$
P_1f(z)=\displaystyle\frac{(-1)^k k}{\pi}\displaystyle\iint_{R^+}\displaystyle\frac{\left[|\zeta|^2-1+\overline{\zeta}(\overline{\zeta}-\overline{z})\right]^{k-1}}{(1-z\overline{\zeta})^{k+1}}\overline{f(\zeta)}\mbox{\rm
d}\xi\mbox{\rm d}\eta$$
and
$$
P_2f(z)=\displaystyle\frac{(-1)^{k+1} k}{\pi}\displaystyle\iint_{R^+}\displaystyle\frac{\left[|\zeta|^2-1+\zeta(\zeta-\overline{z})\right]^{k-1}}{(1-z\zeta)^{k+1}}f(\zeta)\mbox{\rm
d}\xi\mbox{\rm d}\eta.$$
Then, $P_1f$ can be rewritten
$$
P_1f(z)=\displaystyle\frac{(-1)^k k}{\pi}\displaystyle\iint_{\mathbb{D}}\displaystyle\frac{\left[|\zeta|^2-1+\overline{\zeta}(\overline{\zeta}-\overline{z})\right]^{k-1}}
{(1-z\overline{\zeta})^{k+1}}\overline{\widehat{f}(\zeta)}\mbox{\rm
d}\xi\mbox{\rm d}\eta$$
with $\mathbb{D}$ being the unit disc and
$$
\widehat{f}(z)=\left\{\begin{array}{cc}
0,&f\in \mathbb{D}\backslash R^+,\\[3mm]
f(z),&f \in R^+.\end{array}\right.
$$
By \cite{fr1974} and the Theorem 3.1 in \cite{U2007}, we obtain
$$
\parallel P_1f\parallel_{L_p(\mathbb{D})}\leq C(p)\parallel \widehat{f}\parallel_{L_p(\mathbb{D})}=C(p)\parallel f\parallel_{L_p(R^+)},
$$
therefore,
$$\parallel P_1f\parallel_{L_p(R^+)}\leq\parallel P_1f\parallel_{L_p(\mathbb{\mathbb{D}})}\leq C(p)\parallel f\parallel_{L_p(R^+)},
$$
which implies that $P_1f$ is bounded on $L_p(R^+)$ for $p>2$.

Next, since $|1-z\zeta|>|\overline{\zeta}-z|$ and $|1-z\zeta|>|1-z\overline{\zeta}|$ for $z,\zeta\in R^+$, we obtain
$$\begin{array}{ll}
|P_2f(z)|
&\leq C\sum\limits_{m=0}^{k-1}\left(\begin{array}{cc}k-1\\m\end{array}\right)\displaystyle\iint_{R^+}
\displaystyle\frac{(1-|\zeta|^2)^m|\zeta-\overline{z}|^{k-1-m}}{|1-z\zeta|^{k+1}}|f(\zeta)|\mbox{\rm
d}\xi\mbox{\rm d}\eta\\[4mm]
&\leq C\sum\limits_{m=0}^{k-1}\left(\begin{array}{cc}k-1\\m\end{array}\right)\displaystyle\iint_{R^+}
\displaystyle\frac{(1-|\zeta|^2)^m}{|1-z\overline{\zeta}|^{m+2}}|f(\zeta)|\mbox{\rm
d}\xi\mbox{\rm d}\eta\\[4mm]
&=C\sum\limits_{m=0}^{k-1}\left(\begin{array}{cc}k-1\\m\end{array}\right)\displaystyle\iint_{\mathbb{D}}
\displaystyle\frac{(1-|\zeta|^2)^m}{|1-z\overline{\zeta}|^{m+2}}|\widehat{f}(\zeta)|\mbox{\rm
d}\xi\mbox{\rm d}\eta\end{array}$$
with $C$ being a constant. According to the proof of Theorem 3.1 in \cite{U2007}, $P_2$ is bounded on $L_p(\mathbb{D})$ for $p>2$, thus
$$\parallel P_2f\parallel_{L_p(R^+)}\leq\parallel P_2f\parallel_{L_p(\mathbb{D})}\leq
C(p)\parallel \widehat{f}\parallel_{L_p(\mathbb{D})}=C(p)\parallel f\parallel_{L_p(R^+)}.
$$
So, the proof is completed.

\ble For $f\in L_p(R^+)$ and $p>2$,
$$\parallel \Xi_3f\parallel_{L_p(R^+)}\leq C(p)\parallel f\parallel_{L_p(R^+)},$$
where $\Xi_3f$ is defined by \eqref{202111} and $C(p)$ is a constant depending on $p$.
\eele

\p
By \eqref{202111}, suppose  $\Xi_3f(z)=\Upsilon_1f(z)+\Upsilon_2f(z)$ with
$$
\Upsilon_1f(z)=\displaystyle\frac{(-1)^{k} kr^2}{\pi}\displaystyle\iint_{R^+}\displaystyle\frac{\left[|\zeta|^2-r^2+\overline{\zeta}(\overline{\zeta}-\overline{z})\right]^{k-1}}{(r^2-z\overline{\zeta})^{k+1}}\overline{f(\zeta)}\mbox{\rm
d}\xi\mbox{\rm d}\eta$$
and
$$
\Upsilon_2f(z)=\displaystyle\frac{(-1)^{k+1} kr^2}{\pi}\displaystyle\iint_{R^+}\displaystyle\frac{\left[|\zeta|^2 -r^2+\zeta(\zeta-\overline{z})\right]^{k-1}}{(r^2-z\zeta)^{k+1}}f(\zeta)\mbox{\rm
d}\xi\mbox{\rm d}\eta.$$

From the inequality $|\zeta-z|<\frac{1}{r}|r^2-z\overline{\zeta}|$, then
$$\begin{array}{ll}
|\Upsilon_1f(z)|&\leq C_1\sum\limits_{m=0}^{k-1}\left(\begin{array}{cc}k-1\\m\end{array}\right)\displaystyle\iint_{R^+}
\displaystyle\frac{(|\zeta|^2-r^2)^m|\overline{\zeta}-\overline{z}|^{k-1-m}}{|r^2-z\overline{\zeta}|^{k+1}}|f(\zeta)|\mbox{\rm
d}\xi\mbox{\rm d}\eta\\[4mm]
&\leq C_2\sum\limits_{m=0}^{k-1}\left(\begin{array}{cc}k-1\\m\end{array}\right)\displaystyle\iint_{R^+}
\displaystyle\frac{(|\zeta|^2- r^2)^m}{|r^2-z\overline{\zeta}|^{m+2}}|f(\zeta)|\mbox{\rm
d}\xi\mbox{\rm d}\eta.\end{array}$$
Let $\tau=\displaystyle\frac{r}{\zeta}$ and $z=\displaystyle\frac{r}{\mu}$, then $\tau,\,\mu\in R^+$
for $z,\zeta\in R^+$. Also suppose $\tau=\tau_1+i\tau_2$ with $\tau_1,\,\tau_2\in \mathbb{R}$, thus $\mbox{\rm
d}\xi\mbox{\rm d}\eta=\displaystyle\frac{r^2}{|\tau|^4}\mbox{\rm
d}\tau_1\mbox{\rm d}\tau_2$, and
$$\begin{array}{ll}
|\Upsilon_1f(z)|&\leq \displaystyle\frac{C_2|\mu|^{m+2}}{r^2}\sum\limits_{m=0}^{k-1}\left(\begin{array}{cc}k-1\\m\end{array}\right)\displaystyle\iint_{R^+}
\displaystyle\frac{(1-|\tau|^2)^m}{|\tau|^{m+2}|1-\mu\overline{\tau}|^{m+2}}\left|f\left(\displaystyle\frac{r}{\tau}\right)\right|\mbox{\rm
d}\tau_1\mbox{\rm d}\tau_2\\[4mm]
&\leq C_3\sum\limits_{m=0}^{k-1}\left(\begin{array}{cc}k-1\\m\end{array}\right)\displaystyle\iint_{\mathbb{D}}
\displaystyle\frac{(1-|\tau|^2)^m}{|1-\mu\overline{\tau}|^{m+2}}\left|\widetilde{f}(\tau)\right|\mbox{\rm
d}\tau_1\mbox{\rm d}\tau_2\end{array}$$
with
$$\widetilde{f}(\tau)=\left\{\begin{array}{ll}
f\left(\displaystyle\frac{r}{\tau}\right),&\tau\in R^+,\\[4mm]
0,&\tau\in \mathbb{D}\setminus R^+.
\end{array}\right.
$$
Since $f\in L_p(R^+)\Longleftrightarrow \widetilde{f}\in L_p(\mathbb{D})$, and from the above Lemma 4.4,
$$\parallel \Upsilon_1f\parallel_{L_p(R^+)}\leq\parallel \Upsilon_1f\parallel_{L_p(\mathbb{\mathbb{D}})}\leq C(p)\parallel \widetilde{f}\parallel_{L_p(\mathbb{D})}=C(p)\parallel \widetilde{f}\parallel_{L_p(R^+)}\leq C(p)\parallel f\parallel_{L_p(R^+)}.
$$

From $|r^2-z\zeta|>r|\overline{\zeta}-z|$ and $|r^2-z\zeta|>|r^2-z\overline{\zeta}|$, we obtain
$$\begin{array}{ll}
|\Upsilon_2f(z)|
&\leq C_4\sum\limits_{m=0}^{k-1}\left(\begin{array}{cc}k-1\\m\end{array}\right)\displaystyle\iint_{R^+}
\displaystyle\frac{(|\zeta|^2-r^2)^m|\zeta-\overline{z}|^{k-1-m}}{|r^2-z\zeta|^{k+1}}|f(\zeta)|\mbox{\rm
d}\xi\mbox{\rm d}\eta\\[4mm]
&\leq C_5\sum\limits_{m=0}^{k-1}\left(\begin{array}{cc}k-1\\m\end{array}\right)\displaystyle\iint_{R^+}
\displaystyle\frac{(|\zeta|^2-r^2)^m}{|r^2-z\overline{\zeta}|^{m+2}}|f(\zeta)|\mbox{\rm
d}\xi\mbox{\rm d}\eta,\end{array}$$
where $C_i$ for $i=1,\cdots,5$ are constants. According to the above result, $\Upsilon_2$ is bounded on $L_p(R^+)$ for $p>2$, therefore,
$$\parallel \Xi_3f\parallel_{L_p(R^+)}\leq C(p)\parallel f\parallel_{L_p(R^+)}.$$

\ble
Let $f\in L_p(R^+,\mathbb{C})$, $p>2$, $k\in \mathbb{N}$, and
\beq\label{12345}\widetilde{\Pi}_kf\triangleq\displaystyle\frac{\partial^k
\widetilde{T}_k f}{\partial z^k},\eeq then $$\parallel\widetilde{\Pi}_kf\parallel_{L_p(R^+)}\leq C(p)\parallel f\parallel_{L_p(R^+)},$$
where $C(p)$ is a constant depending on $p$.
\eele

\p
From \eqref{6.3}, \eqref{2021999} and \eqref{20212021} for $l=k-1$,  and differentiating with respect to $z$, we obtain
\beq\label{2021888}
\displaystyle\frac{\partial^k
\widetilde{T}_k[f](z)}{\partial z^k}
=\Xi_1f(z)+\Xi_2f(z)+\Xi_3f(z)+\Xi_4f(z)+\Xi_5f(z),
\eeq
where $\Xi_1,$ $\Xi_2$, $\Xi_3$ are given by \eqref{202109},\eqref{202110} and \eqref{202111}, respectively. Moreover,
$$
\Xi_4f(z)=\displaystyle\frac{1}{\pi}\displaystyle\iint_{R^+}\Delta(z,\zeta)f(\zeta)
\mbox{\rm
d}\xi\mbox{\rm d}\eta,
$$
$$
\Xi_5f(z)=\displaystyle\frac{1}{\pi}\displaystyle\iint_{R^+}\Delta(z,\overline{\zeta})\overline{f(\zeta)}
\mbox{\rm
d}\xi\mbox{\rm d}\eta
$$
with
$$\begin{array}{ll}
\Delta(z,\zeta)=&\displaystyle\sum\limits_{j=0}^{k-2}\displaystyle\frac{(-1)^{k-j-1}(k-1)!}{j!(k-j-2)!}
\left(\zeta-z+\overline{\zeta-z}\right)^{k-j-2}h^*(z,\zeta)\\[5mm]
&+\displaystyle\sum\limits_{j=0}^{k-1}\left(\begin{array}{cc}k-1\\j\end{array}\right)(-1)^{k-j}
\left(\zeta-z+\overline{\zeta-z}\right)^{k-j-1}\partial_zh^*(z,\zeta),
\end{array}
$$
where $h^*(z,\zeta)=g^*(z,\zeta)$ for $l=k-1$ given by \eqref{20212021}. Here we know $h^*(z,\zeta)$ and $\partial_z h^*(z,\zeta)$ are bounded for $z,\zeta\in R^+$, which means that $\Delta(z,\zeta)$ is bounded for $z,\zeta\in R^+$. Thus,
$$
\parallel \Xi_4f\parallel_{L_p(R^+)}\leq C(p)\parallel f\parallel_{L_p(R^+)}.
$$

Similarly,
$$
\parallel \Xi_5f\parallel_{L_p(R^+)}\leq C(p)\parallel f\parallel_{L_p(R^+)}.
$$
In terms of Lemma 4.3-4.5 and \eqref{2021888}, the proof is completed.

 With the above preliminaries, we can solve the following Schwarz problem for a generalized inhomogeneous polyanalytic equation on basis of the properties of the operators $\widetilde{T}_k$.

{\bf  Schwarz Problem}\ Find $w\in W^{p,n}(R^+)$ as a solution to the complex differential equation as follows
\beq\label{9.1}\begin{array}{ll}
\displaystyle\frac{\partial^nw}{\partial\overline{z}^n}&+\sum\limits_{j=1}^n q_{1j}(z)\displaystyle\frac{\partial^nw}{\partial\overline{z}^{n-j}\partial z^j}+\sum\limits_{j=1}^n q_{2j}(z)\displaystyle\frac{\partial^n\overline{w}}{\partial z^{n-j}\partial \overline{z}^j}\\[4mm]
&+\displaystyle\sum\limits_{l=0}^{n-1}\sum\limits_{m=0}^l\left[ a_{ml}(z)\displaystyle\frac{\partial^lw}{\partial\overline{z}^{l-m}\partial z^m}+
b_{ml}(z)\displaystyle\frac{\partial^l\overline{w}}{\partial z^{l-m}\partial \overline{z}^m}\right]=f(z),\ \ z\in R^+,
\end{array}\eeq
with $a_{ml}$, $b_{ml}$, $f\in L_p(R^+,\mathbb{C})$, $p>2$ and $q_{1j}$, $q_{2j}$ are measurable bounded functions
 \beq\label{9.2}
 \displaystyle\sum\limits_{j=1}^n|q_{1j}(z)|+|q_{2j}(z)|<q_0<1,
 \eeq
 which satisfying the boundary conditions
\beq\label{9.3}
\mbox{Re}\displaystyle\frac{\partial^lw}{\partial\overline{z}^l}=\gamma_l,\ \ z\in \partial R^+,\ \ l=0,1,\cdots,n-1.
\eeq

Similarly in \cite{ay2018}, from \eqref{55.4}, \eqref{6.1}, Lemma 3.2 and Theorem 3.1, we have the following result.

\bth
The above Schwarz problem is equivalent to the following equation
\beq\label{9.7}
(I+\widetilde{\Xi}+\widetilde{K})g=\chi(z),
\eeq
with
$w=\widetilde{T}_{n}g+\Phi(z)$,
$$
\widetilde{\Xi}g=\sum\limits_{j=1}^n\left(q_{1j}\widetilde{\Pi}_jg+q_{2j}\overline{\widetilde{\Pi}_jg}\right),
$$
$$
\widetilde{K}g=\displaystyle\sum\limits_{l=0}^{n-1}\sum\limits_{m=0}^l\left[ a_{ml}(z)\displaystyle\frac{\partial^m\widetilde{T}_{n+m-l}g}{\partial z^m}+
b_{ml}(z)\displaystyle\frac{\partial^m\overline{\widetilde{T}_{n-l+m}g}}{\partial \overline{z}^m}\right],$$
$$\begin{array}{ll}
\Phi(z)=&i\sum\limits_{k=0}^{n-1}(z+\overline{z})^kc_k+\displaystyle\sum\limits_{k=0}^{n-1}\frac{(-1)^k}{k!}\bigg\{\frac{1}{2\pi
i}\int_{\Gamma_1\bigcup\Gamma_2}\gamma_k(\zeta)\left(\zeta-z+\overline{\zeta-z}\right)^kH_1(z,\zeta)\displaystyle\frac{\mbox{\rm
d}\zeta}{\zeta}\\[7mm]
&\ \ \ \ +\displaystyle\frac{1}{\pi
i}\int_{[\varpi,w]\cup[r,1]}\gamma_k(\zeta)\left(\zeta-z+\overline{\zeta-z}\right)^kH_2(z,\zeta)\mbox{\rm
d}\zeta\bigg\},
\end{array}
$$
and
$$\begin{array}{ll}
\chi(z)=f-&\sum\limits_{j=1}^n \left(q_{1j}(z)\displaystyle\frac{\partial^n\Phi}{\partial\overline{z}^{n-j}\partial z^j}+ q_{2j}(z)\displaystyle\frac{\partial^k\overline{\Phi}}{\partial z^{n-j}\partial \overline{z}^j}\right)\\[4mm]
&+\displaystyle\sum\limits_{l=0}^{n-1}\sum\limits_{m=0}^l\left[ a_{ml}(z)\displaystyle\frac{\partial^l\Phi}{\partial\overline{z}^{l-m}\partial z^m}-
b_{ml}(z)\displaystyle\frac{\partial^l\overline{\Phi}}{\partial z^{l-m}\partial \overline{z}^m}\right].
\end{array}$$
Here, $\widetilde{T_k}$, $\widetilde{\Pi}_k$ are defined by \eqref{6.1} and \eqref{12345}, respectively. $H_1$, $H_2$  are given by \eqref{5.1} and \eqref{5.2} for $\alpha=1$, respectively.
\eeth

\bth
Suppose the condition \beq\label{789}q_0\max\limits_{1\leq j\leq n}\parallel\widetilde{\Pi}_j  \parallel_{L_p(R^+)}<1,\eeq
then the equation \eqref{9.1} with the conditions \eqref{9.2} and \eqref{9.3} has a solution of $w=\widetilde{T}_n g+\Phi(z)$, where $g\in L_p(R^+)$, $p>2$
is a solution of the singular integral equation \eqref{9.7}.
\eeth

\p
By Lemmas 4.1-4.2 and Arzela Ascoli theorem, it is easy to know that $\widetilde{K}$ is compact. Moreover, Lemma 4.6 and \eqref{789} imply that $I+\widetilde{\Xi}$ is an invertible operator on $L^p(R^+)$. Then $I+\widetilde{\Xi}+\widetilde{K}$ is of Fredholm type with index zero. So the singular integral equation \eqref{9.7} is solvable. Thus, by Theorem 4.1, the Schwarz problem \eqref{9.1}$-$\eqref{9.3} is solvable and its solution is $w=\widetilde{T}_n g+\Phi(z)$, where $g$ is the solution of \eqref{9.7} given by $\left(I+\widetilde{\Xi}+\widetilde{K}\right)^{-1}\chi(z)$. The proof is completed.

\section*{Acknowledgements}
This work was supported by National Natural Science
Foundation of China (11801570).

\end{document}